# Workflow for High-Fidelity Dynamic Analysis of Structures with Pile Foundation

Amin Pakzad
*Civil and Environmental Engineering, PhD candidate, University of Washington, Seattle. United States, amnp95@uw.edu*
**Pedro Arduino**
*Civil and Environmental Engineering, Professor, University of Washington, Seattle, United States, parduino@uw.edu*
Wenyang Zhang
*Research Assoc., Texas Advanced Computing Center, University of Texas at Austin*, wzhang@tacc.utexas.edu
Ertugrul Taciroglu
*Civil and Environmental Engineering, Professor, University of California Los Angeles, etacir@ucla.edu*

ABSTRACT: The demand for high-fidelity numerical simulations in soil-structure interaction analysis is on the rise, yet a standardized workflow to guide the creation of such simulations remains elusive. This paper aims to bridge this gap by presenting a step-by-step guideline proposing a workflow for dynamic analysis of structures with pile foundations. The proposed workflow encompasses instructions on how to use Domain Reduction Method for loading, Perfectly Matched Layer elements for wave absorption, soil-structure interaction modeling using Embedded interface elements, and domain decomposition for efficient use of processing units. Through a series of numerical simulations, we showcase the practical application of this workflow. Our results reveal the efficacy of the Domain Reduction Method in reducing simulation size without compromising model fidelity, show the precision of Perfectly Matched Layer elements in modeling infinite domains, highlight the efficiency of Embedded Interface elements in establishing connections between structures and the soil domain, and demonstrate the overall effectiveness of the proposed workflow in conducting high-fidelity simulations. While our study focuses on simplified geometries and loading scenarios, it serves as a foundational framework for future research endeavors aimed at exploring more intricate structural configurations and dynamic loading conditions.

## 1 INTRODUCTION

The demand for high-fidelity numerical simulations in soil structure interaction analysis is increasing. Finite element simulation is widely used for conducting such simulations. However, while previous research efforts in literature have focused on specific problems, there is currently no standard workflow to guide the creation of these simulations for common structures.

This paper aims to address this gap by proposing a step-by-step guideline for high fidelity dynamic analysis of structures with pile foundations. This guideline provides comprehensive directions covering several aspects, including how to apply loads, select the domain of interest, absorb reflecting waves, model soil-structure interaction, and utilize domain decomposition for effective utilization of multiple processing units. At the end, a potential application of this workflow is demonstrated by providing results of simulations.

## 2 WORKFLOW

In this work, a step-by-step workflow is proposed for creating high-fidelity simulations for dynamic analysis of structures with pile foundations. This is achieved by breaking the process into smaller parts.

### 2.1 *Domain Reduction Method*

Near-fault regions exhibit unique ground motions, making simulating earthquake motions at a regional scale essential. Although using high-performance computing accelerates these regional-scale models, they remain computationally heavy. Therefore, there is a need to reduce computational complexity. The Domain Reduction Method (DRM) facilitates the integration of earthquake sources, considering soil-structure interaction and accommodating complex 3D wavefields in a computationally efficient manner. DRM is a two-step numerical approach for modeling earthquake-induced ground motion in three dimensions. Initially developed by Bielak et al. (2003) and Yoshimura et al. (2003), DRM simplifies complex scenarios by dividing them into two stages.

First, an auxiliary problem is solved, simulating a large-scale domain encompassing the seismic wave field and the vicinity of the structure under examination. This stage, known as the background problem, removes localized features to focus on the broader seismic wave field.

In the second step, the domain is reduced to the immediate surroundings of the structure, loaded with forces equivalent to the seismic excitation computed in the first step. This enables detailed analysis of the structure's response, including nonlinear effects. Figure 1 shows the schematic representation of the domain reduction method steps.

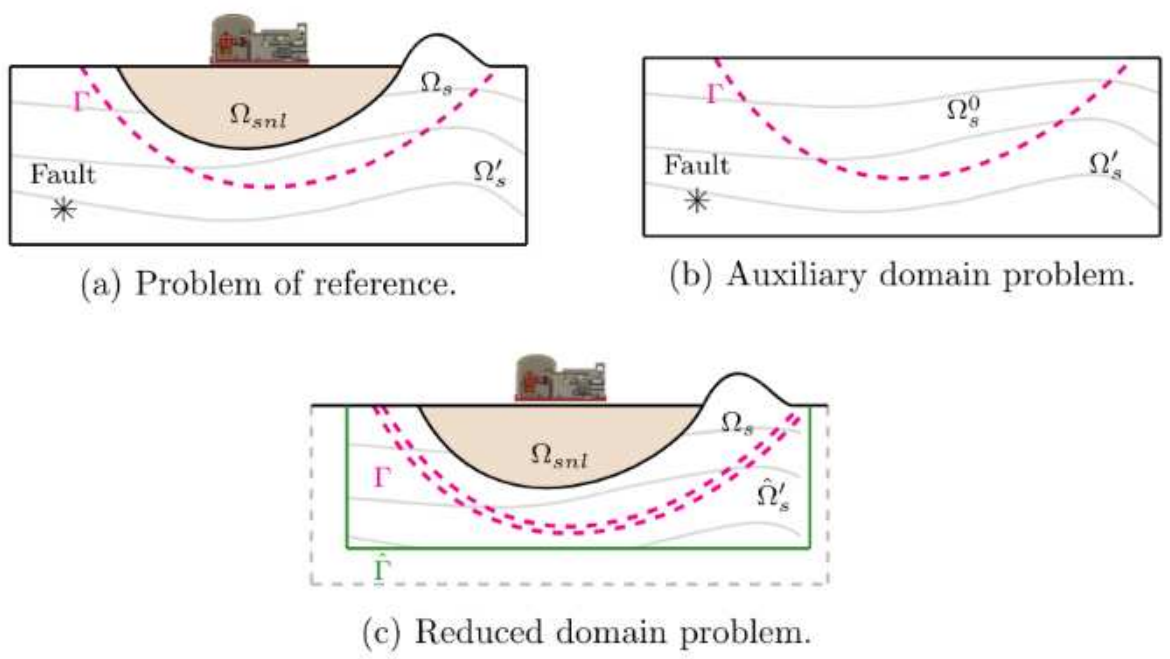

Figure 1 Schematic representation of the domain reduction method steps Korres et.al (2023).

## 2.2 *Perfectly Matched Layer*

The utilization of a finite element method necessitates transforming the infinite physical domain into a finite one for computational purposes. This process involves establishing absorbing boundary conditions at the domain edges. In these simulations, an unsplit-field Perfectly-Matched-Layer (PML) element, as described in Fathi et al. (2015), Wang (2019) and Trono (2023) is utilized.

The fundamental concept of the Perfectly Matched Layer (PML) is to construct a layer surrounding the computational domain's edges to absorb all outgoing waves. This layer is designed to match the properties of the surrounding material, ensuring that waves passing through it are smoothly absorbed without reflecting.

The formulation of the PML is based on the concept of complex coordinate stretching, where the usual physical coordinates of the element are substituted with complex-valued coordinates. This innovative approach significantly reduces the computational workload involved in dynamic earthquake analysis by eliminating the need to extend the computational domain just to absorb waves. Figure 2 shows the schematic of the PML absorbing layer in the simulations.

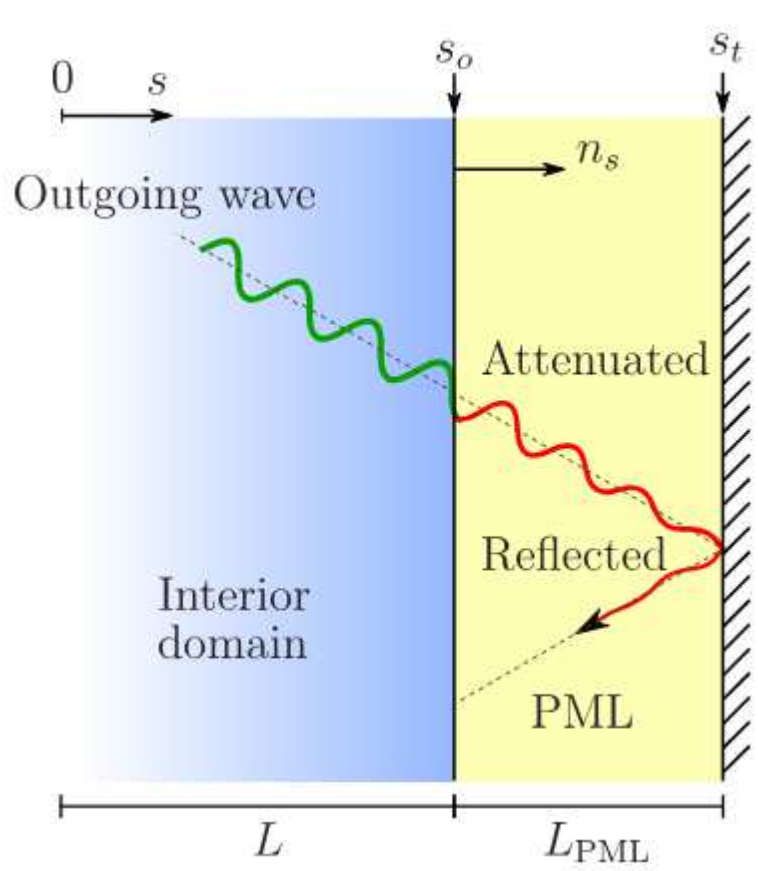

Figure 2 Schematic of the PML absorbing Element Fathi et.al (2015).

## 2.3 *Embedded Interface Elements*

The analysis of lateral loading on piles encompasses various techniques, from simple two-dimensional beam over elastic foundation theory to highly refined three-dimensional (3D) finite element (FE) models. Accurately defining boundary conditions at the soil-pile interface is crucial for capturing realistic pile behavior. However, determining suitable interface conditions, including interface friction, contact stiffness, and soil-pile relative displacements, can be challenging, especially in 3D models due to meshing difficulties. Incorporating a simple pile into the model can alter the geometry and necessitate mesh refinement, thereby increasing computational demands.

The concept of embedded beam elements, formulated for the elastic analysis of piles in FE models, offers a solution to reduce computational demands. In this approach, the displacement field of the pile is expressed in terms of the displacement interpolation of the surrounding solid (soil) elements. Sadek and Shahrour (2004) developed an embedded element for piles, assuming perfect bonding between a beam embedded in its surrounding solid. In their formulation, the bonding condition is kinematically enforced at the beam nodes, assuming that the displacement field along the beam can be defined in terms of the solid displacement field expressed by standard interpolation functions. However, this approach may lead to non-convergent solutions, particularly when the beam and solid nodes are located close to each other.

To address this issue, Turello et al. (2016) explicitly defined an interaction surface along which the kinematic condition of compatibility is enforced. This is achieved using the principle of virtual works in a weak sense, ensuring that the relative displacements between the beam and the solid result in zero virtual work for any admissible system of virtual interaction forces. Turello et al. (2017) expanded this formulation to consider elasto-plastic behavior at the interface between the beam and solid element. The embedded element proposed by Turello et al. is suitable for modeling soil-pile systems, providing accuracy and ease of use in numerical modeling, especially for dynamic analysis of such problems. In this workflow, the idea is to utilize embedded interface elements for soil-pile interaction, as proposed by Turello et al. For this purpose, the work by Ghofrani (2018) is used here.

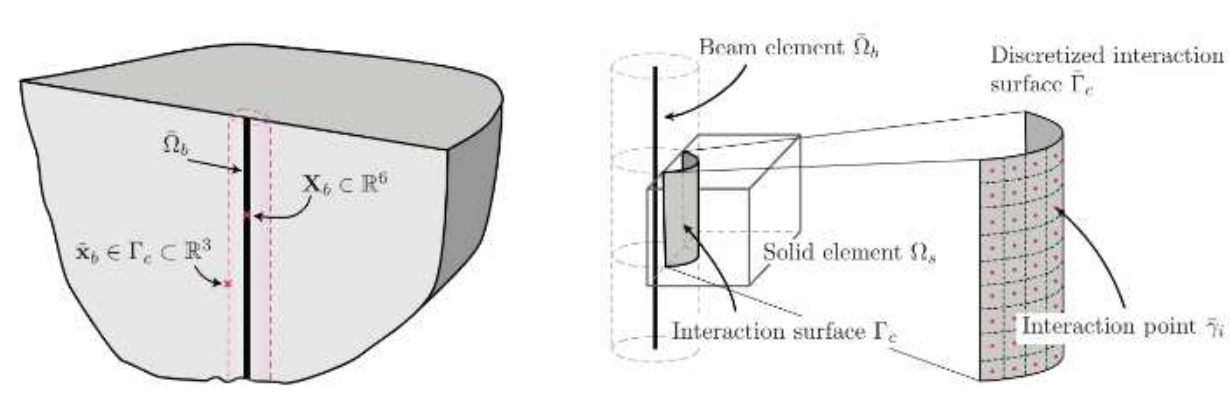

Figure 3 Schematic of a 1D pile foundation buried in a 3D soil domain with an explicit interaction surface Ghofrani (2018).

## 2.4 *Domain Decomposition*

Using DRM, PML, and Embedded Interface elements reduces computational demand, but the scale of models required for creating a representative high-fidelity simulation remains significant, often taking several days for a single run. To address this problem the primary objective of the workflow is to model complex systems feasibly, necessitating a reduction in runtime through optimal resource utilization. Employing parallel computing addresses this challenge by distributing the analysis across multiple CPUs using domain decomposition within the workflow.

However, employing domain decomposition and parallel computing introduces additional challenges regarding how to perform the decomposition effectively. The optimal approach to decomposition needs careful consideration. As discussed, our model domain is divided into three primary categories: the Regular domain, which encompasses the soil and structure of interest; the DRM domain, consisting of a single layer of elements surrounding the Regular domain, utilized for force generation; and the PML layer, typically comprising 2 to 6 layers of elements surrounding the DRM layer, designed to absorb outgoing waves reflecting from the soil or structure within the Regular domain. In 3D each hexahedral PML element within the PML layer imposes a computational demand greater than that of Regular domain elements. This is because PML formulations require at least 12 degrees of freedom per node. Conversely, the stresses and displacements on the DRM layer elements are not conducive to straightforward post-processing, as they undergo significant deformation and stress under loads. Embedded elements also impose constraints requiring interaction points and corresponding soil elements to be associated with one partition.

Consequently, for domain decomposition, the whole domain must be partitioned into three/four main subdomains (Regular, DRM, and PML, pile if present), with each domain parsed into different processors. The figure illustrates a sample model mesh partitioned into different processors, with each color representing a distinct processor.

Figure 4 shows a typical mesh employed in generating high-fidelity simulations utilizing the DRM method, PML elements, and structural connections to the soil domain with embedded elements. In Figure 5, the same mesh is displayed partitioned into distinct processors, where each color corresponds to a different processor.

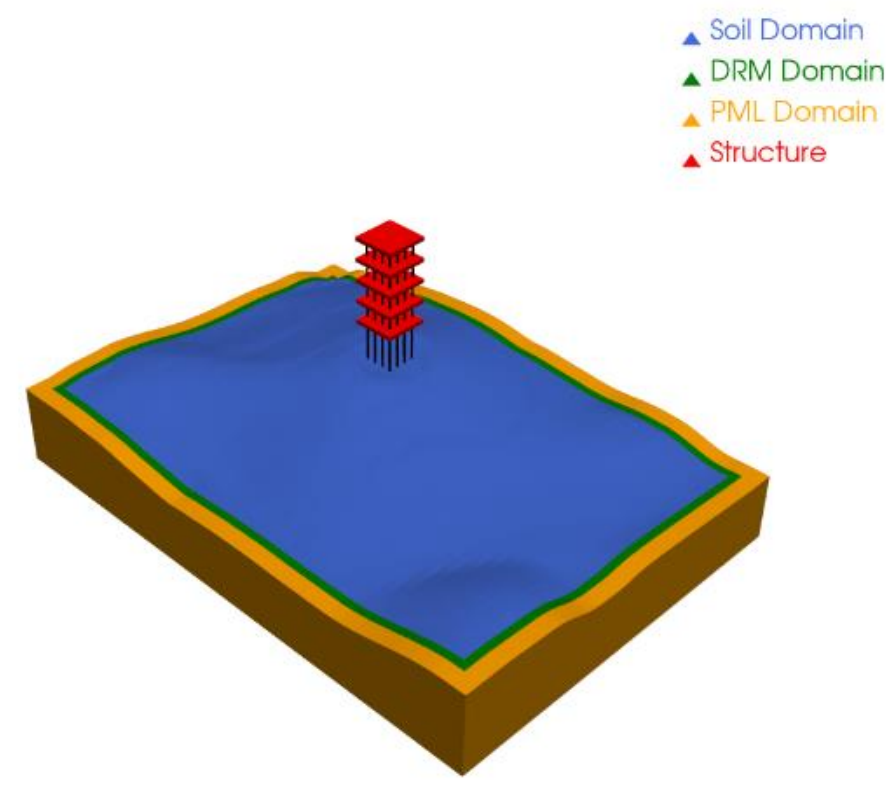

Figure 4 Typical mesh used in the proposed workflow including all the domains.

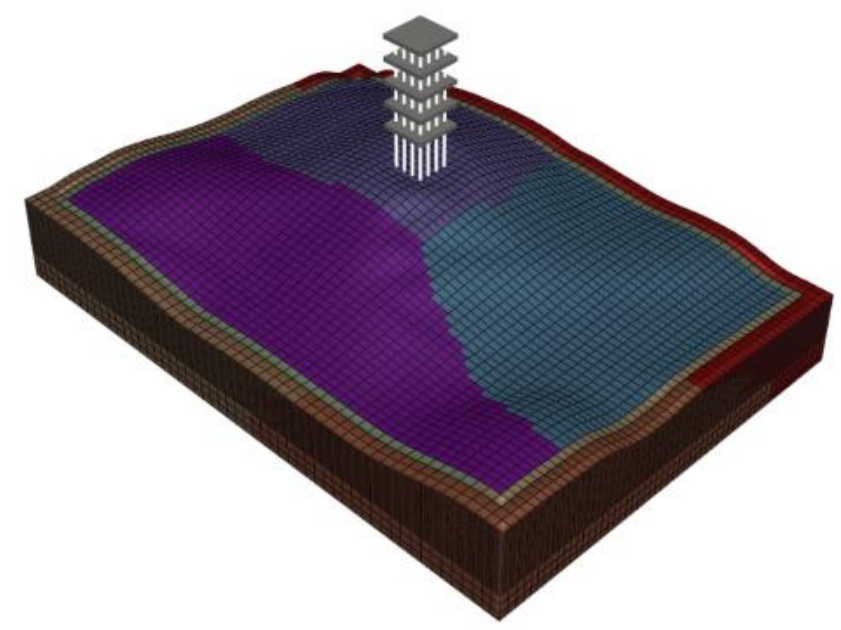

Figure 5 Example of partitioned mesh for distribution among multiple processors.

## 3 NUMERICAL SIMULATIONS

Presented below are a couple of numerical simulations and their results to show the workflow's application. The finite element software OpenSees (Open System for Earthquake Engineering Simulations McKenna et.al (2010)) is utilized for performing all simulations. For DRM loading, an existing load pattern in OpenSees is employed, which reads HDF5 format to load accelerations and velocity time histories on nodes in the DRM layer. The forces on these nodes are then computed automatically using the stiffness matrix of the elements. Regarding PML elements, the core Fortran code used for implementing these elements in Abaqus by Zhange et al. (2019) has been adapted for implementation in OpenSees. To verify the workflow, five sets of numerical simulations are considered.

To generate a Domain ReductionMethod (DRM) loading, one must determine the acceleration and displacement values at nodes positioned along the boundaries of the region of interest. In this study, we adopt a straightforward approach by examining a soil profile where waves originate from a specific depth and propagate vertically towards the surface. The computation of DRM loading entails utilizing the accelerations and displacements of nodes distributed throughout the domain. For the simplified 1D wave propagation scenario, we utilize an elastic material characterized by an elastic modulus of 200 MPa and a Poisson's ratio of 0.3. These elastic properties, along with the DRM load derived from the 1D analysis, serve as the basis for all subsequent numerical simulations. A Ricker wavelet is employed as the source load, applied at the base of the model. Figure 6 shows the time history of the ricker wavelet at the base of the model and it it's frequency content.

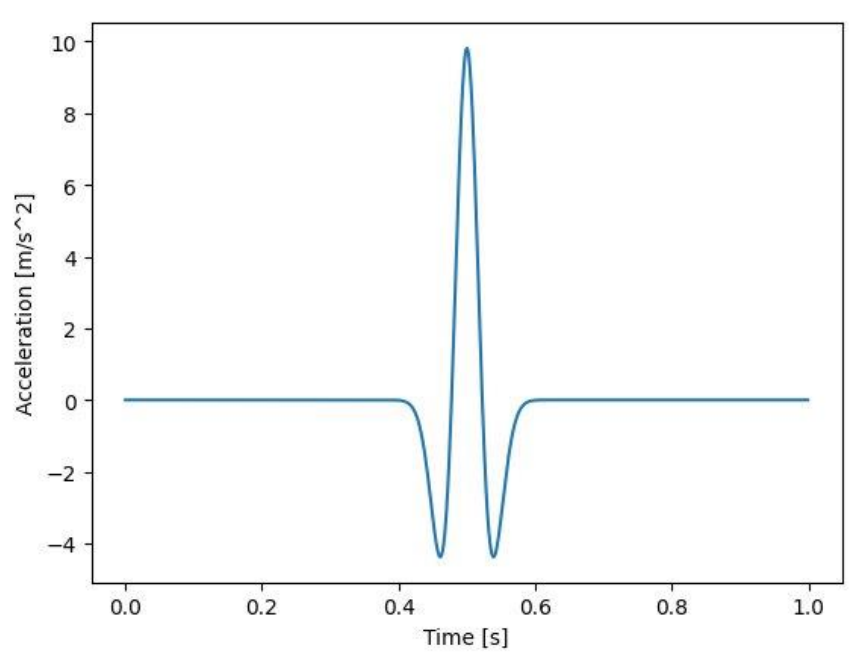

a) Time history

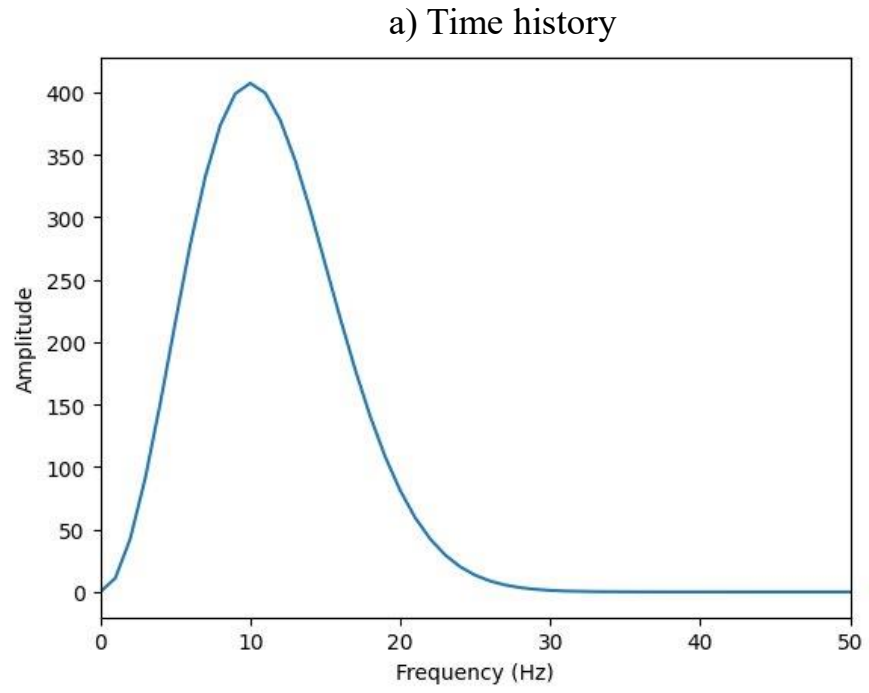

b) Frequency content

Figure 6 Time history and frequency content of the loading used for 1D wave propagation problem.

Five sets of simulations are conducted to investigate the efficacy of the workflow. For each set, a numerical model is developed. The analysis of this model is carried out in two ways. In the first simulation, after applying the DRM layer for loading, the boundaries of the model are fixed. In the second simulation, two layers of PML are added to the sides to absorb waves.

### 3.1 *Level Ground Simulation*

As a first set, we consider a 3D flat model with dimensions of 16 meters width, 40 meters depth, and 1 meter length. The DRM loading for this model is obtained from a 1D site response analysis using a single soil column. Recorded accelerations and displacements within the soil column are employed to create the DRM loading. Figure 7 illustrates the mesh used for this model, with a mesh size of 0.5 meters. The objective of this numerical simulation is to assess the effectiveness of DRM loading in modeling wave propagation within the domain, as well as to evaluate the performance of the PML elements.

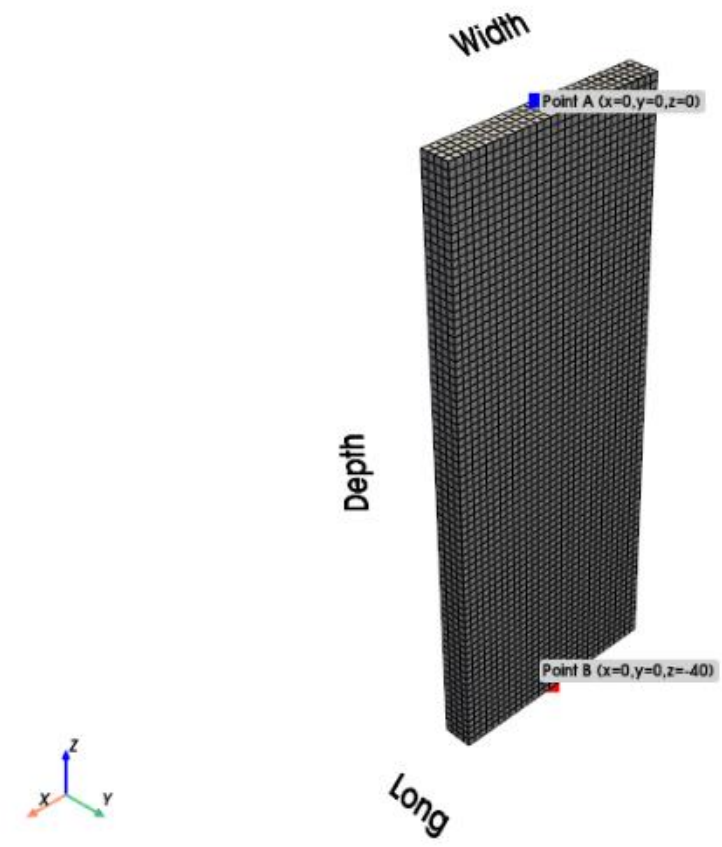

Figure 7 Mesh used for level ground simulation and the selected control points.

Two different nodes are selected at different depths—one at the surface (point A) and the other at the point where the wave enters the domain (point B)—for comparison. Figure 8 shows the comparison of accelerations at points A and B. As the wave enters from point B into the domain, traverses the surface, and reflects to the bottom of the elastic layer, the maximum acceleration amplitude at point A is observed to be twice that at point B. This result aligns with the expectations due to the free boundary conditions at the surface, consistent with what was recorded in the 1-dimensional model. In the second approach, two layers of PML elements are added after the DRM layer to the sides and bottom of the layer, and the model is analyzed again. The response of accelerations at points A and B remains unchanged compared to the model with fixed boundaries. This lack of change is expected because there is no structure at the free surface to receive waves and reflect waves with a different frequency. Therefore, the addition of PML elements is expected to have no impact on the model, as observed.

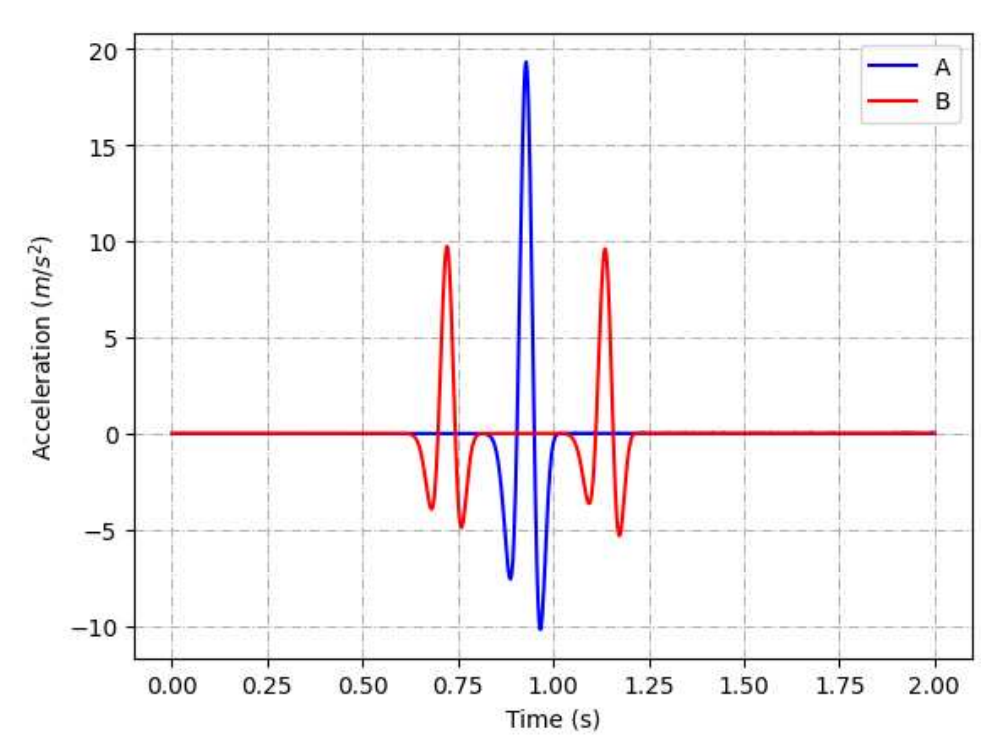

Figure 8 Comparison of acceleration on the selected control points A and B.

### 3.2 *Reduced Domain Level Ground Simulation*

One of the primary advantages of the domain reduction method (DRM) lies in its ability to shrink the domain of analysis while preserving wave propagation and other properties of the waves being simulated. The concept revolves around narrowing the scope of analysis until nonlinear elements within the domain become significant, influencing the surrounding areas. To validate the efficacy of DRM in a reduced domain, a secondary simulation was performed using a distinct 3D model with altered dimensions derived from the initial simulation, featuring fewer elements. In this model, the width remained consistent at 16 meters, similar to the first simulation, while the depth decreased to 8 meters and the length increased to 6 meters. Once more, simulations were conducted employing fixed boundaries and utilizing perfectly matched layer (PML) elements at the boundaries. This downsized model served as the foundational framework for subsequent simulations involving structural elements within the domain. Figure 9 shows the mesh with the control point A utilized in this model. A typical simulation used 3 partitions for the PML layer, 1 for DRM and a 3 partitions for the soil domain.

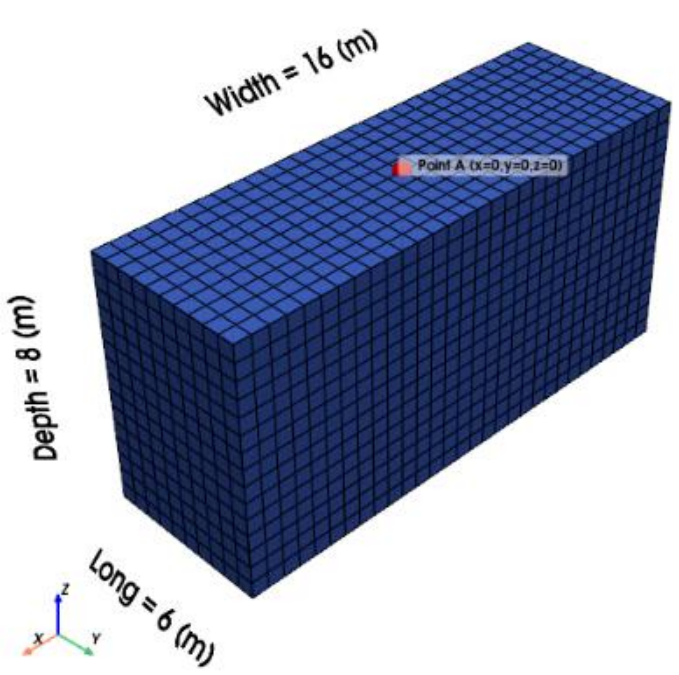

Figure 9 Mesh used for reduced level ground simulation and the selected control points.

Figure 10 displays the acceleration time history for control point A. The acceleration time history matches exactly with that obtained at the surface for point A in the initial simulation. Once more, utilizing fixed boundaries versus PML boundaries makes no difference in the model, as there is no new structural interactions from the surface of the model to the boundaries.

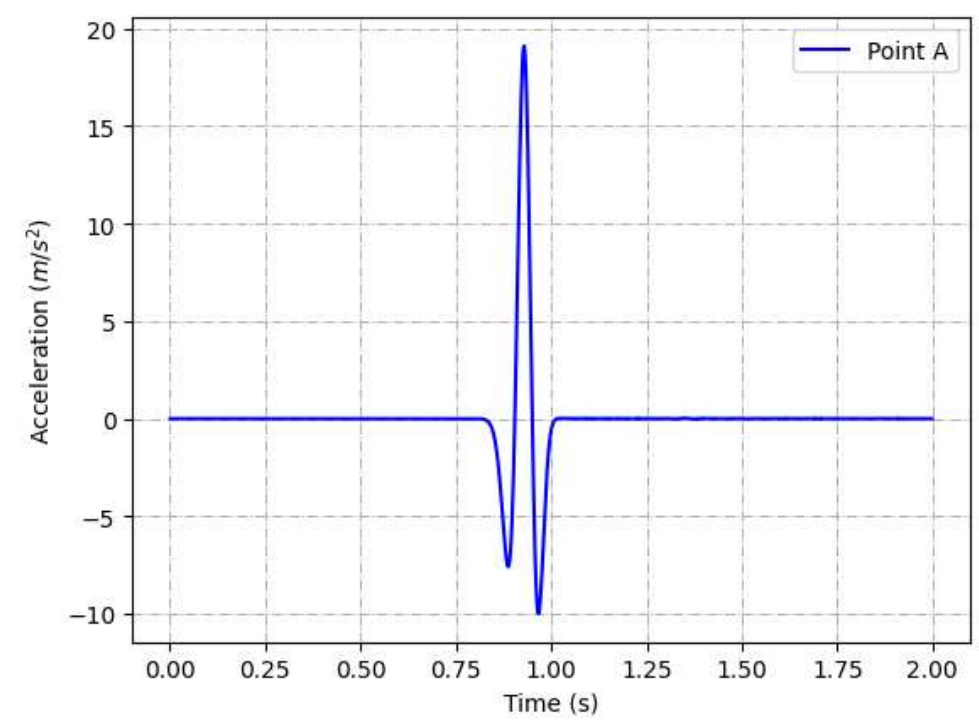

Figure 10 Acceleration time history of the Point A at the surface of the reduce level ground simulation.

### 3.3 *Lumped Mass on Pile Foundation*

For the third set of simulations, a pile foundation with a lumped mass on top is incorporated. The presence of this pile alters the original wave propagation in the domain. Specifically, the wave traverses through the pile, causing it to oscillate and reflect new waves in the soil domain. To create this model, a reduced flat model from the second simulation is utilized as the base model. For the pile foundation, a pile with a 5-meter embedded length and a height of two meters above the ground, and a lumped mass equivalent to 50 tons on the top, considered. The pile is simulated using a beam-column element with the same stiffness as a circular section pile with a diameter of 1.0 meter, an elastic modulus of 10 GPa, and a Poisson's ratio of 0.3.

To connect the pile to the soil mesh, embedded interface elements was utilized, assuming perfect bonding between the soil and the pile. Using this element enabled us to connect the soil mesh to the pile without any modifications on the soil mesh. Figure 11 shows the mesh used in this simulation for both the soil and pile parts.

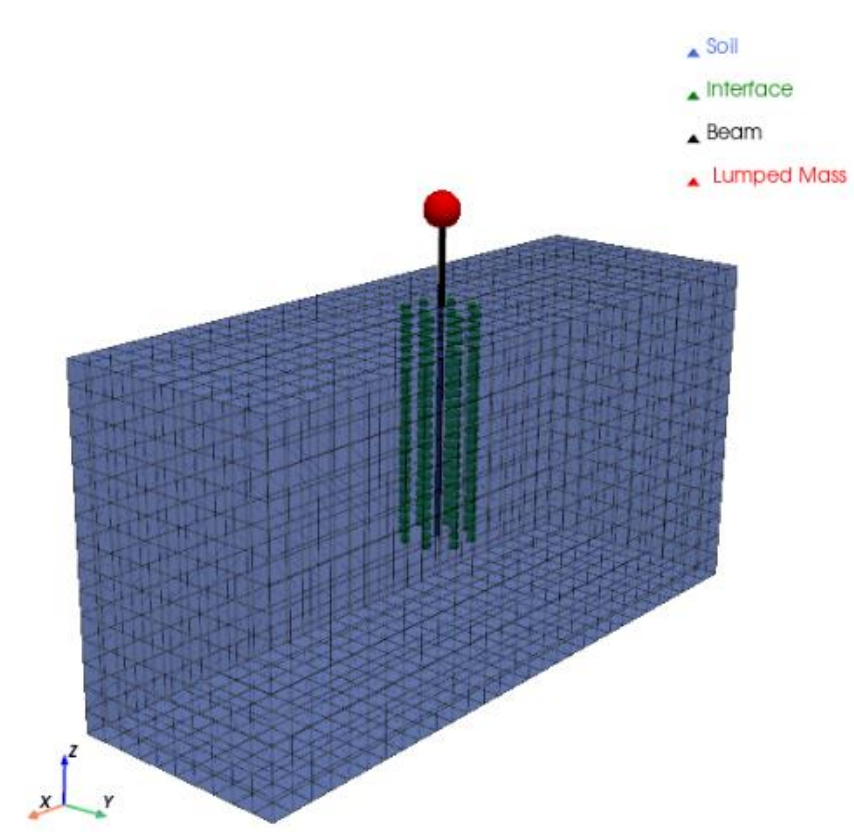

Figure 11 Mesh used in the third simulation set simulation for both the soil and pile parts.

Figure 12 shows the displacement time history of the top of the pile. Comparing the results for the cases with fixed boundaries and PML boundaries, it is observed that the amplitude of the pile oscillation decreases in the case with PML boundaries after the wave passes through the pile and produces additional oscillations. However, for the fixed boundaries case, the amplitude remains unchanged, and since there is no damping in the numerical model, the pile continues to oscillate for the duration of the simulation.

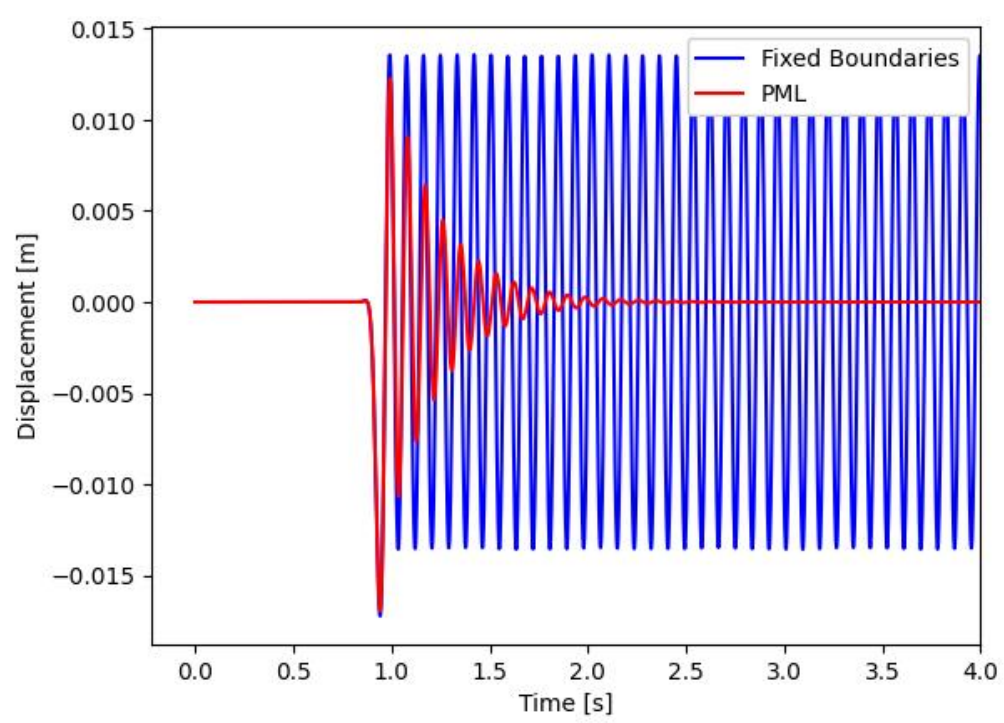

Figure 12 Displacement time history of the pile head in the third simulation.

Figure 13 and Figure 14 illustrate contours of acceleration on a vertical section of the domain for simulations with fixed and PML boundaries respectively at time = 1.5 (s). Comparing the acceleration contours and considering the decrease in amplitude for the simulation with PML boundaries, it can be observed that the reason for the amplitude decrease with PML boundaries is that using PML boundaries allows waves to pass through the boundaries. This simulates infinite boundary conditions in the model, altering all boundary conditions, which aligns with a model that has an infinite domain from the sides and domain.

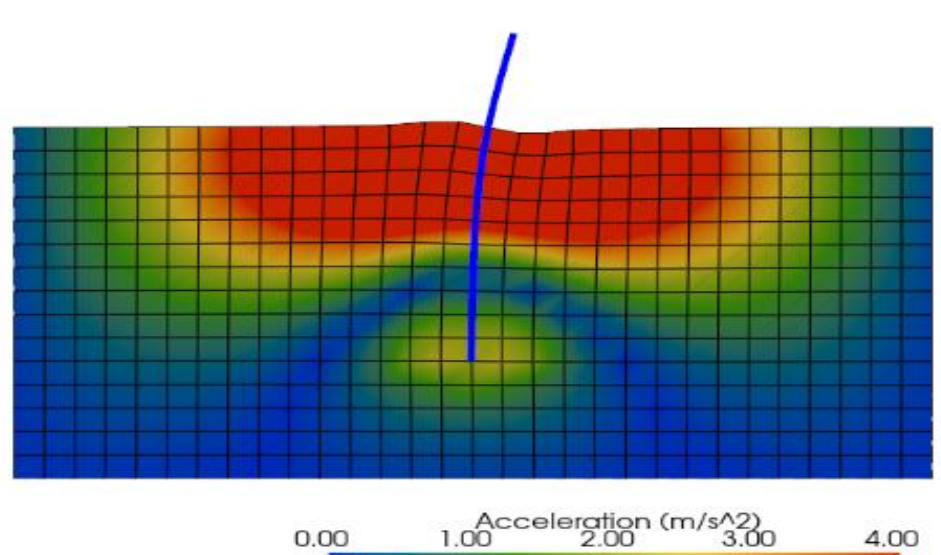

Figure 13 Contours of acceleration magnitude on a vertical section of the domain for simulation with fixed boundaries at time=1.5 s.

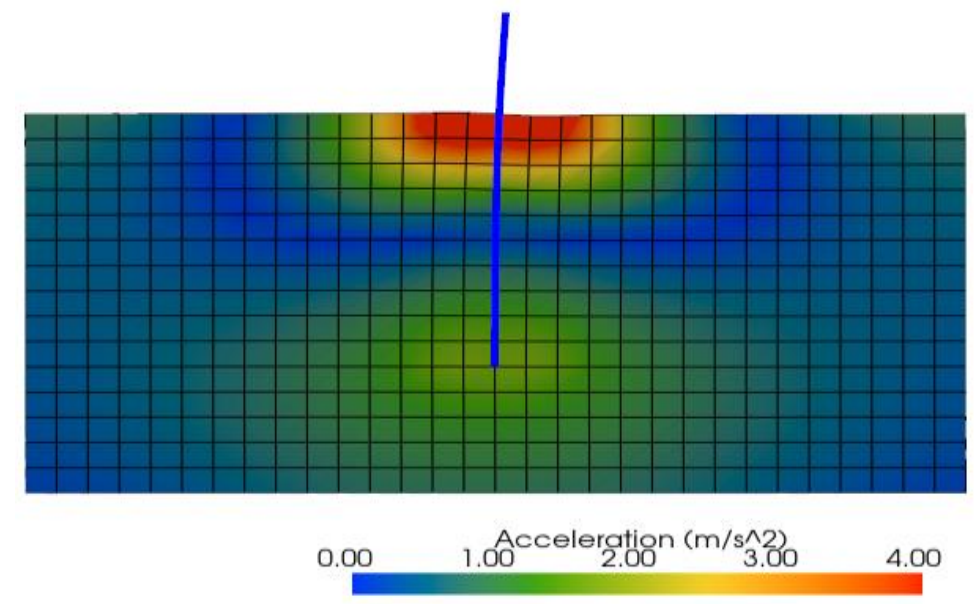

Figure 14 Contours of acceleration magnitude on a vertical section of the domain for simulation with PML boundaries at time=1.5 s.

The results of this simulation using DRM loading and PML elements demonstrate the applicability of the workflow in modeling the dynamic behavior of pile foundations.

### 3.4 *Concrete Structure Simulation*

For the fourth model a concrete structure with width of 3 meters and length of 3 meter and height of 3 meters inserted at the center of the base model from the second simulation attached to the soil layer at the surface. For creating this structure a Elastic material with Elastic modulus of 25 GPa and Poisson ratio of 0.22 and density of 2400 $kg/m^3$ is used as representative of concrete. Figure 15 shows the mesh used in this simulation.

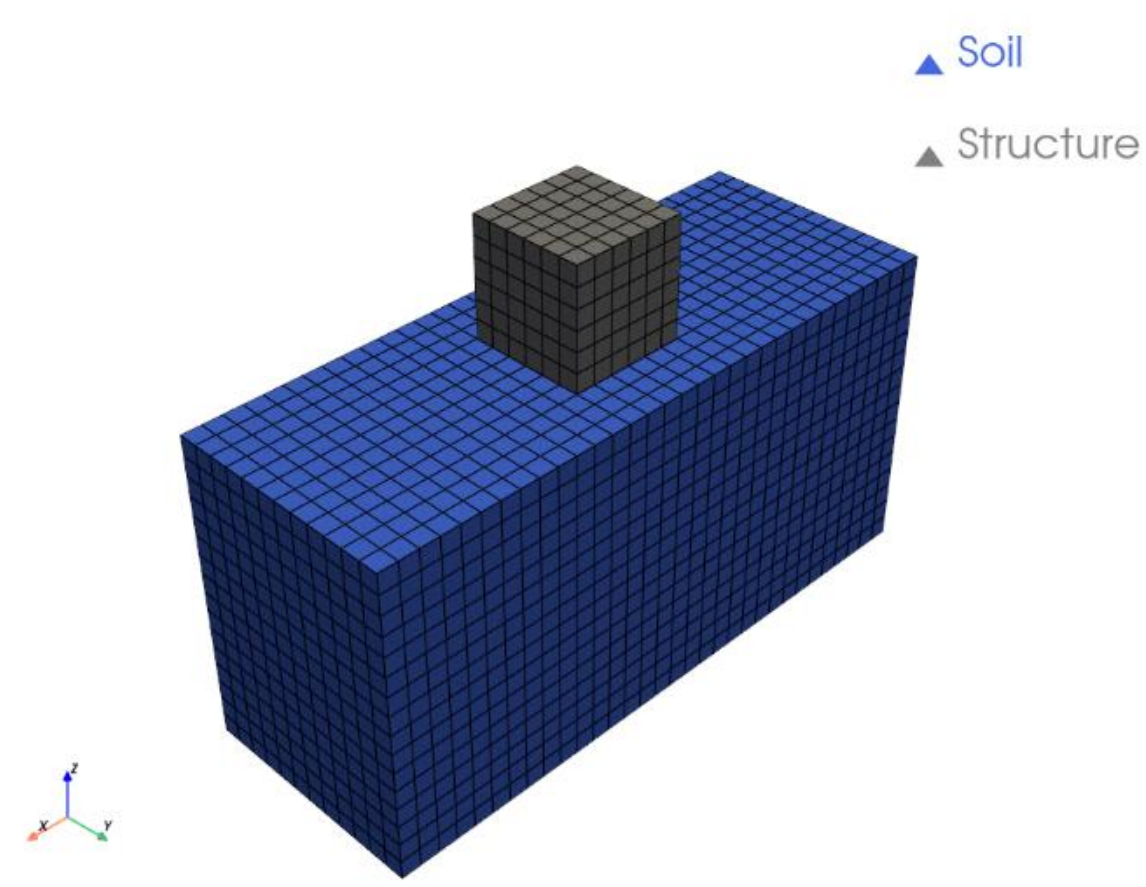

Figure 15 Mesh used in the fourth simulation showing a concrete structure on the soil layer.

Figures 11 and 12 depict the displacement time history of the top and bottom of the structure's center. The structure begins to oscillate due to the loading, and, similar to the third simulation set, the amplitudes of the simulations using the PML boundaries decrease to zero, indicating infinite domain boundary conditions. Conversely, the amplitudes of the simulations with fixed boundaries remain constant throughout the analysis.

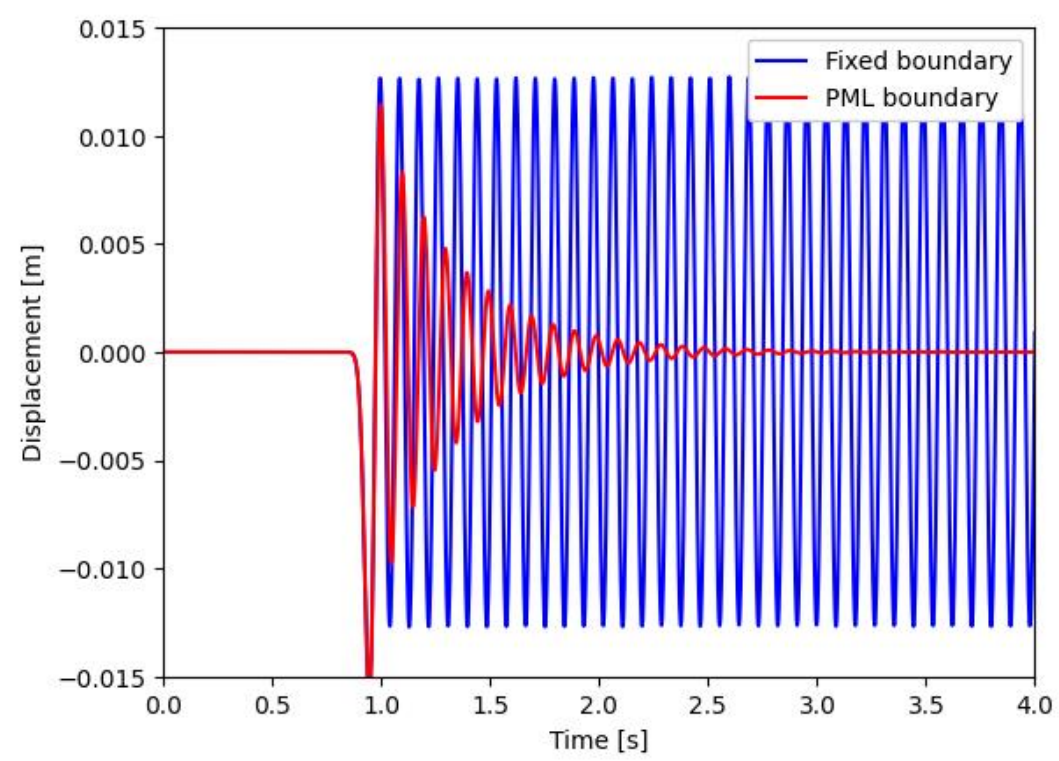

Figure 16 Displacement time history of the top of structure.

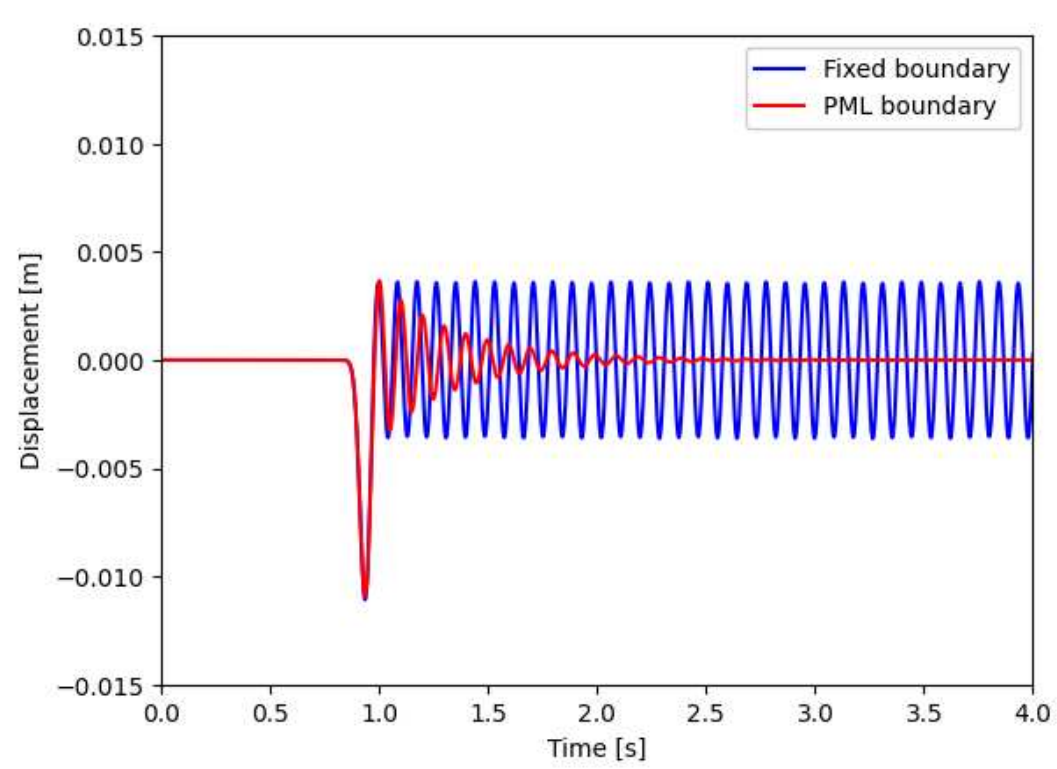

Figure 17 Displacement time history of the bottom of structure.

## 3.5 *Multiple Degrees of Freedom Oscillator on a Foundation*

For this simulation set, we use a beam model of the structure, which mirrors the configuration of the previous simulation. This beam incorporates multiple degrees of freedom and ultimately rests on an isolated foundation placed on level ground. To construct these models, a square foundation measuring 3 meters in width, 3 meters in length, and 0.5 meters in depth is positioned within the reduced soil domain from the second simulation. For the foundation, an elastic material with an elastic modulus of 25 GPa, a Poisson's ratio of 0.22, and a density of 2400 $kg/m^3$ is utilized as a representation of concrete material.

The beam in this simulation has a height of 2.75 meters and lateral stiffness equivalent to the lateral stiffness of the building's cross-section. The equivalent mass of the structure is distributed among the nodes of the beam. To connect the beam to the foundation, an embedded interface element is employed, embedding a length of 0.25 meters of the beam into the foundation. Figure 18 shows the soil mesh, beam, foundation, interface elements and the control points selected for this simulation.

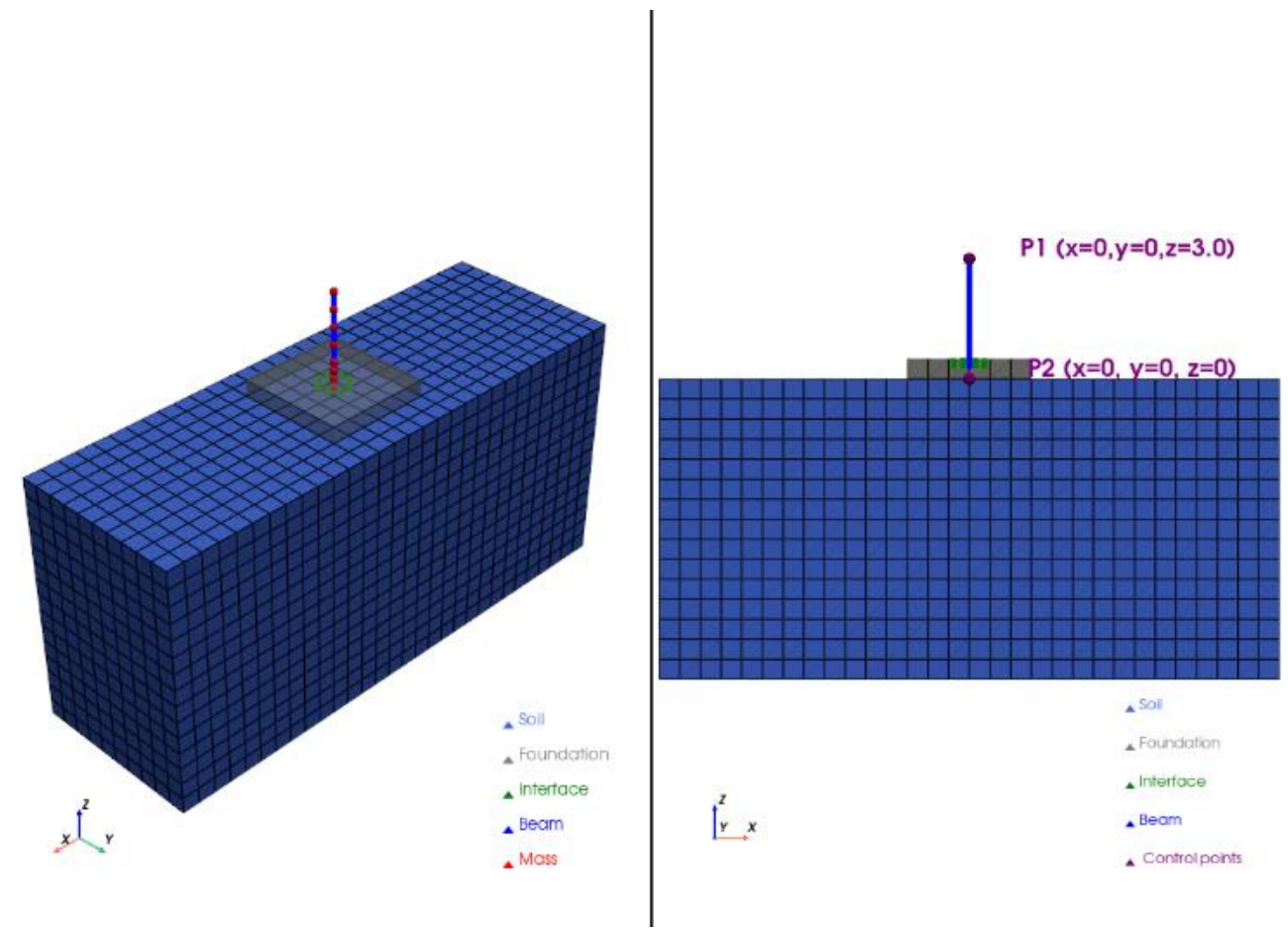

Figure 18 Mesh used for the fifth simulation set.

Figure 19 and Figure 20 display the displacement time histories of control points P1 and P2, representing the top and bottom of the structure, respectively, from the previous simulation set. The displacement time histories exhibit similar behavior for the multiple degrees of freedom oscillator when compared to the equivalent concrete structure from which this model is derived, both in simulations with fixed and PML boundary conditions.

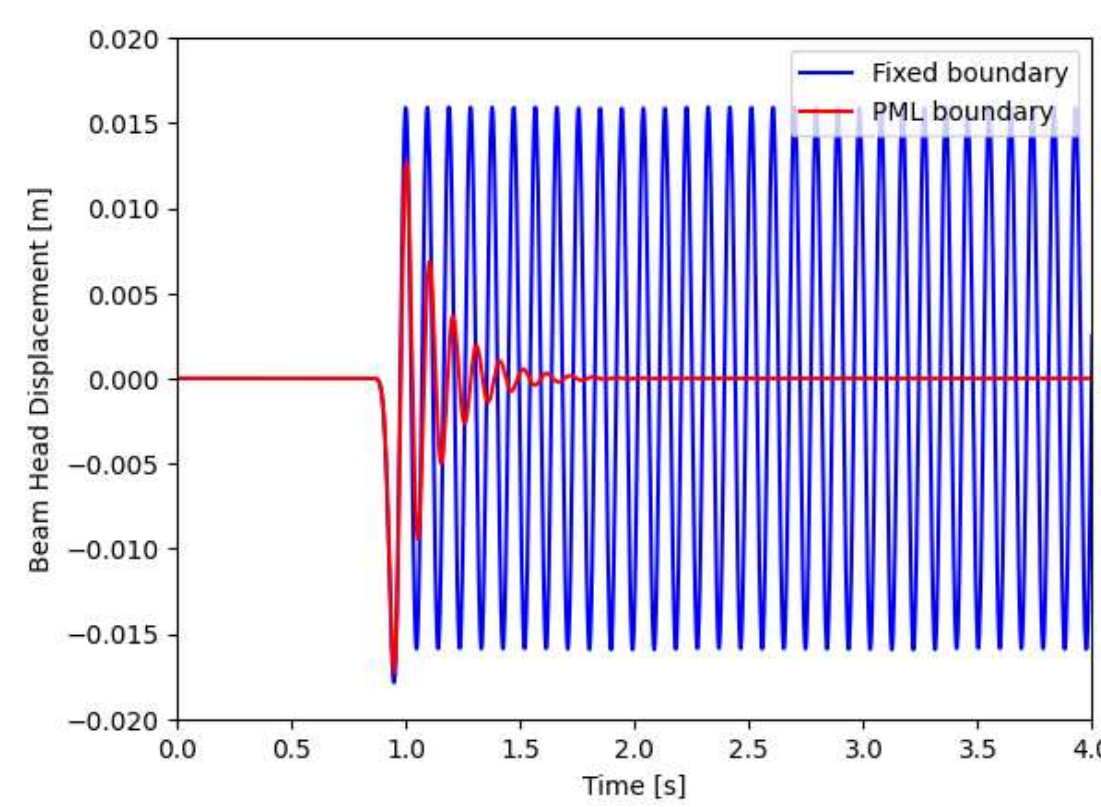

Figure 19 Displacement time histories of control points P1.

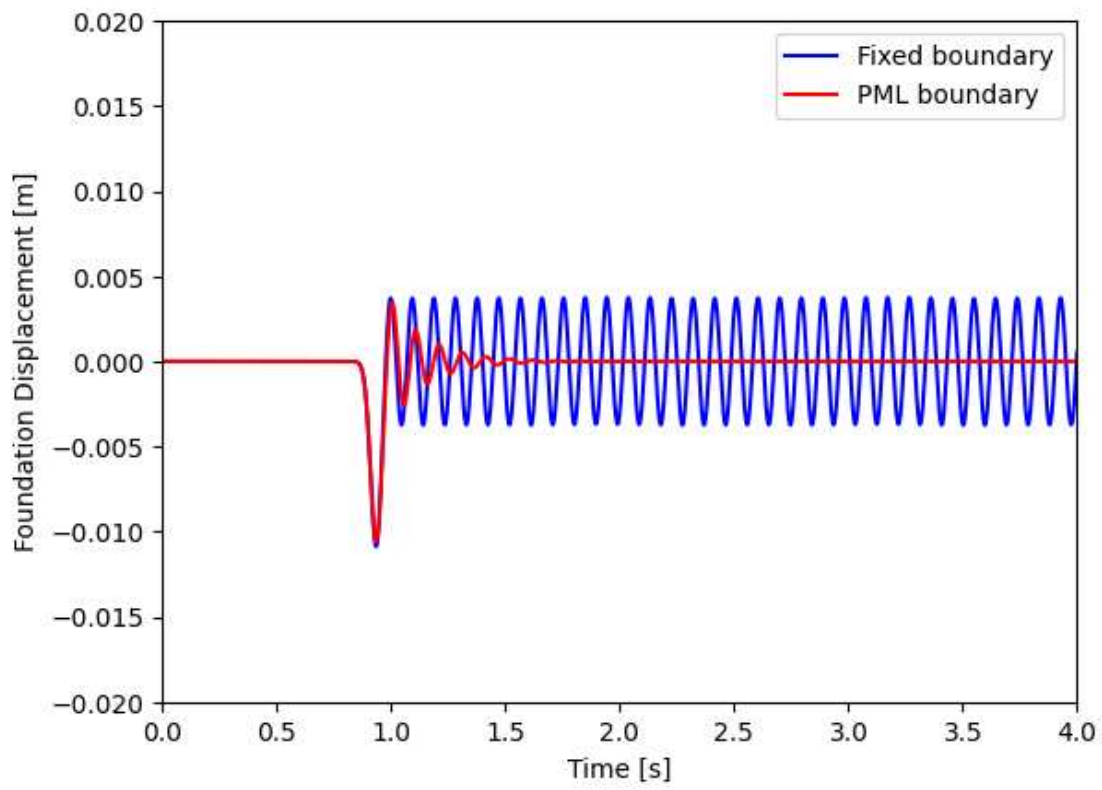

Figure 20 Displacement time histories of control points P2.

All simulations demonstrate the potential efficacy of utilizing the DRM method, PML elements, and Embedded Interface elements in creating high-fidelity numerical simulations for dynamic structural analysis across various aspects. This workflow could be further utilized for modeling more complex structures.

## 4 CONCLUSIONS

In this study, we propose a workflow for conducting high-fidelity dynamic analysis of structures with pile foundations. This workflow integrates various components, including the Domain Reduction Method, Perfectly Matched Layer elements, Embedded Interface elements, and Domain Decomposition. The efficacy of this workflow was investigated through a series of numerical simulations following the proposed methodology.

The outcomes of our simulations revealed several findings:

- The Domain Reduction Method effectively reduced the simulation size while preserving the essential features of the model.
- Perfectly Matched Layer elements successfully modeled the infinite domain in the simulation, enhancing its accuracy.
- Embedded Interface elements efficiently connected different structures to the soil domain, facilitating simulations.
- The overall workflow showed effectiveness in conducting complex high-fidelity numerical simulations.

While this study showed the application of the proposed workflow, it is important to acknowledge its limitations. The scope of our simulations was limited to simplified structural geometries, an elastic soil domain and loading scenarios. Future research will explore more complex structural configurations and dynamic loading conditions to further validate our findings.

In conclusion, this research has provided insights into the capabilities of high-fidelity dynamic analysis of structures and has set the stage for future investigations in this field.

## 5 ACKNOWLEDGEMENTS

The authors of this paper extend their gratitude to DesignSafe, PEER, and SimCenter for providing the funding support for this research. The contributions from these organizations will play a crucial role in the successful execution of this research project.